\pdfoutput=1
\RequirePackage{ifpdf}
\ifpdf 
\documentclass[pdftex]{sigma}
\else
\documentclass{sigma}
\fi

\numberwithin{equation}{section}

\newtheorem{Theorem}{Theorem}[section]
\newtheorem{Proposition}[Theorem]{Proposition}
\newtheorem{Lemma}[Theorem]{Lemma}

\theoremstyle{definition}
\newtheorem{Remark}[Theorem]{Remark}
\newtheorem{Definition}[Theorem]{Definition}

\newcommand{\ta}{\theta}
\newcommand{\C}{\mathbb C}
\newcommand{\E}{\mathbb E}
\newcommand{\N}{\mathbb N}
\newcommand{\Z}{\mathbb Z}

\begin{document}

\newcommand{\arXivNumber}{2307.12921}

\renewcommand{\thefootnote}{}

\renewcommand{\PaperNumber}{052}

\FirstPageHeading

\ShortArticleName{An elliptic algebra and elliptic multinomial theorem}

\ArticleName{An Algebra of Elliptic Commuting Variables\\ and an Elliptic Extension of the Multinomial Theorem\footnote{This paper is a~contribution to the Special Issue on Basic Hypergeometric Series Associated with Root Systems and Applications in honor of Stephen C.~Milne's 75th birthday. The~full collection is available at \href{https://www.emis.de/journals/SIGMA/Milne.html}{https://www.emis.de/journals/SIGMA/Milne.html}}}

\Author{Michael J. SCHLOSSER}

\AuthorNameForHeading{M.J.~Schlosser}

\Address{Fakult\"at f\"ur Mathematik, Universit\"at Wien,\\ Oskar-Morgenstern-Platz~1, A-1090 Vienna, Austria}
\Email{\href{mailto:michael.schlosser@univie.ac.at}{michael.schlosser@univie.ac.at}}
\URLaddress{\url{http://www.mat.univie.ac.at/~schlosse/}}

\ArticleDates{Received December 26, 2024, in final form June 23, 2025; Published online July 06, 2025}

\Abstract{We introduce an algebra of elliptic commuting variables involving a base~$q$, nome~$p$, and $2r$ noncommuting variables. This algebra, which for $r=1$
reduces to an algebra considered earlier by the author, is an elliptic extension of the well-known algebra of $r$ $q$-commuting variables. We present a multinomial theorem valid as an identity in this algebra, hereby extending the author's previously obtained elliptic binomial theorem to higher rank. Two essential ingredients are a consistency relation satisfied by the elliptic weights and the Weierstra{\ss} type $\mathsf A$ elliptic partial fraction decomposition. From the elliptic multinomial theorem we obtain, by convolution, an identity equivalent to Rosengren's type~$\mathsf A$ extension of the Frenkel--Turaev ${}_{10}V_9$ summation. Interpreted in terms of a weighted counting of lattice paths in the integer lattice $\mathbb Z^r$, this derivation of Rosengren's $\mathsf A_r$ Frenkel--Turaev summation constitutes the first combinatorial proof of that fundamental identity.}


\Classification{05A10; 11B65; 33D67; 33D80; 33E90}

\renewcommand{\thefootnote}{\arabic{footnote}}
\setcounter{footnote}{0}

\section{Introduction}\label{sec:intro}
A standard and powerful tool in algebraic combinatorics
is the identification of a class of combinatorial objects with
words of noncommuting variables in some monoid.
Such a correspondence is convenient when
the combinatorics of the objects translates, on the algebraic
side, to commutation relations satisfied by the variables.
For instance, semi-standard Young tableaux are in
one-to-one correspondence with words in the \textit{plactic monoid}
(cf.\ \cite{LS}), the monoid consisting of all words
in an alphabet of totally ordered variables modulo
\textit{Knuth equivalence}. Another example, one which we
shall consider in this paper in the ``elliptic'' setting, is the
correspondence between positively oriented lattice paths
of length $m$ in the $r$-dimensional integer lattice and words
of length $m$ in $r$ variables built from the successive
steps of the path. In the context of weighted enumeration
(typically with respect to the \textit{area} of the path,
or a similar statistic),
commutation relations can be used in connection with
\textit{normalization} to determine the weight of a path.
A~well-known example---intimately tied to
$q$-combinatorics---is the monoid of $q$-commuting
variables in~${X_1,\dots,X_r}$, satisfying the
$q$-commutation relations $X_jX_i=qX_iX_j$, for
$1\le i<j\le r$, where~$q$ is an indeterminate.
(In Section~\ref{sec:ellalg}, we consider an~elliptic
algebra that generalizes exactly this monoid.)
Rather than restricting ourselves just to a monoid,
as we are interested in generating functions
we consider formal linear combinations of
elements of the monoid, thus we work in a~unital
associative algebra. 
This has the advantage that by exploiting certain
properties (such as associativity) of the elements
in the algebra one is able to prove 
combinatorial identities
by entirely algebraic means.
For instance, working in the just mentioned
algebra of $q$-commuting variables,
a proof of the $q$-multinomial theorem
(which is a particular summation theorem for
multivariate basic hypergeometric series)
can be given easily.

The aim of this paper is to develop appropriate algebraic machinery
that fits with the theory of $\mathsf A_r$ elliptic hypergeometric series.
Specifically, we introduce an algebra of elliptic commuting variables
involving a base $q$, nome $p$, and $2r$ noncommuting
variables. This algebra, which for~${r=1}$
reduces to an algebra considered earlier by the author~\cite{S1},
is an elliptic extension of the aforementioned well-known
algebra of $r$ $q$-commuting variables. We present
a multinomial theorem valid as an identity in this algebra,
hereby extending the author's previously obtained elliptic
binomial theorem from \cite{S1} to higher rank.
From the elliptic multinomial theorem, we obtain, by convolution,
an identity equivalent to Rosengren's type $\mathsf A$ extension
of the Frenkel--Turaev ${}_{10}V_9$ summation.
Interpreted in terms of a weighted counting of lattice paths
in the integer lattice $\mathbb Z^r$, this derivation of Rosengren's
$\mathsf A_r$ Frenkel--Turaev summation constitutes
the first combinatorial proof of that fundamental identity.

After having explained the gross outline of the paper, we
briefly explain what elliptic hypergeometric series are,
as the main application of the elliptic algebra introduced
in this paper concerns identities for such series.
\textit{Elliptic hypergeometric series} form a
natural extension of
ordinary and of basic hypergeometric series.
Consider the series
$S=\sum_{k\ge 0}c(k)$, and~${g(k):=c(k+1)/c(k)}$
being the quotient of two consecutive terms of $S$.
The series $S$ is by definition an ordinary (or ``rational'')
hypergeometric series if the ratio $g(k)$ is a rational
function in the summation index $k$. Similarly, $S$ is a
basic (or ``$q$-'', or ``trigonometric'') hypergeometric
series if $g(k)$ is a rational function in $q^k$ (the base
$q$ usually satisfying $|q|<1$). Finally, $S$ is an elliptic
hypergeometric series if $g(k)$ is an \textit{elliptic}
function in $k$ (by which one understands a
complex-valued function that is doubly-periodic and
meromorphic.) Elliptic hypergeometric series made
their first implicit appearance in 1987 in the work of
the mathematical physicists Date, Jimbo, Kuniba, Miwa
and Okado~\cite{DJKMO} as elliptic $6$-$j$ symbols,
representing elliptic solutions of the Yang--Baxter equation.
Ten years later, Frenkel and Turaev~\cite{FT}, by exploiting the
tetrahedral symmetries of those $6$-$j$ symbols and making the
expressions explicit, wrote out the first identities for (what they called)
``modular hypergeometric series'' (now commonly called
\textit{elliptic hypergeometric series}). In particular, they discovered
what is now called the ${}_{12}V_{11}$
transformation (an elliptic extension of Bailey's very-well-poised
${}_{10}\phi_9$ transformation) and, by applying specialization,
the ${}_{10}V_9$ summation (which is an elliptic extension of
Jackson's very-well-poised~${}_8\phi_7$ summation).

We start with explaining some important notions from the
theory of elliptic hypergeometric series (cf.\ \cite[Chapter~11]{GR}
and \cite{R21}) which we shall need.
Let $\C^\times:=\C\setminus\{0\}$. Let the \textit{modified
Jacobi theta function} (in short: \textit{theta function}) with
argument $x$ and fixed nome $p$ be defined by~${
\ta(x)=\ta (x; p):= (x; p)_\infty (p/x; p)_\infty}$, $
\ta (x_1, \dots, x_m): = \prod^m_{i=1} \ta (x_i)$,
where $ x, x_1, \dots, x_m, p\in \C^\times$, ${|p| < 1}$, and
$(x; p)_\infty=\prod^\infty_{k=0}\bigl(1-x p^k\bigr)$ is an infinite
$p$-shifted factorial.

The theta function satisfies the following simple properties,
namely the \textit{inversion}
\begin{subequations}
\begin{gather}\label{eq:inv}
\ta(x)=-x \ta(1/x),
\end{gather}
the \textit{quasi-periodicity}
\begin{gather}\label{p1id}
\ta(px)=-\frac 1x \ta(x),
\end{gather}
\end{subequations}
and the \textit{three-term addition formula}
(cf.\ \cite[p.~451, Example 5]{WW})
\begin{gather}\label{addf}
\ta(xy,x/y,uv,u/v)-\ta(xv,x/v,uy,u/y)
=\frac uy \ta(yv,y/v,xu,x/u).
\end{gather}
The addition formula in \eqref{addf} is a special case of the
following more general identity due to Weierstra{\ss}
(cf.\ \cite[p.~451, Example 3]{WW}),
which we refer to as the \textit{elliptic partial fraction identity}
of type $\mathsf A$:
let $a_1,\dots,a_r,b_1,\dots,b_r\in\C^\times$, then
\begin{gather}\label{pfdA}
\sum_{i=1}^r\frac{\prod_{j=1}^r\ta(a_i/b_j)}{\prod_{j\neq i}\ta(a_i/a_j)}=0,
\end{gather}
under the assumption that the \textit{elliptic balancing condition}
$a_1\cdots a_r=b_1\cdots b_r$ holds.
The addition formula in~\eqref{addf} is a rewriting of the $r=3$
special case of \eqref{pfdA}. While the relation in~\eqref{addf}
serves as key ingredient in the theory of elliptic hypergeometric
series, the partial fraction decomposition in \eqref{pfdA}
is underlying the theory of multivariate elliptic hypergeometric
series associated to the root system $\mathsf A_{r-1}$
(cf.\ \cite{R04,RW}). Indeed, in the theory of (multivariate)
elliptic hypergeometric series inductive proofs and functional
equations typically make use of the identities in \eqref{addf}
and \eqref{pfdA}
(or, in the setting of root systems other than $\mathsf A_{r-1}$,
of other suitable elliptic partial fraction identities which exist).

Now define the {\em theta shifted factorial} (or
{\em $q,p$-shifted factorial}) by
\begin{gather*}
(a;q,p)_n := \begin{cases}
\displaystyle \prod^{n-1}_{k=0} \ta \bigl(aq^k\bigr),& n = 1, 2, \dots,\cr
1,& n = 0,\cr
\displaystyle 1/\prod^{-n-1}_{k=0} \ta \bigl(aq^{n+k}\bigr),& n = -1, -2, \dots.
\end{cases}
\end{gather*}
For compact notation, we write
\begin{gather*}
(a_1, a_2, \dots, a_m;q, p)_n := \prod^m_{k=1} (a_k;q,p)_n.
\end{gather*}
Notice that for $p=0$ one has $\ta (x;0) = 1-x$,
in which case $(a;q, 0)_n = (a;q)_n$
is a {\em $q$-shifted factorial} in base $q$ (cf.\ \cite{GR}).

Notice that
\[
(pa;q,p)_n=(-1)^na^{-n}q^{-\binom n2} (a;q,p)_n,
\]
which follows from repeated use of \eqref{p1id}.
A list of other useful identities for manipulating the
$q,p$-shifted factorials is given in \cite[Section~11.2]{GR}.

By definition, a function $g(u)$ is {\em elliptic} if it is
a doubly-periodic meromorphic function of the complex variable $u$.

Without loss of generality, one may assume
(see \cite[Theorem~1.3.3]{R21}) that
\begin{gather*}
g(u)=\frac{\ta(a_1q^u,a_2q^u,\dots,a_sq^u)}
{\ta(b_1q^u,b_2q^u,\dots,b_sq^u)} z
\end{gather*}
(i.e., $g$ is an abelian function of some degree $s$,
cf.\ \cite{Ba}),
for a constant $z$ and some
$a_1,a_2,\dots,a_s$, $b_1,b_2,\dots,b_s$, $q,p\in\C^\times$ with $|p|<1$,
where the elliptic balancing condition, namely
\[
a_1a_2\cdots a_s=b_1b_2\cdots b_s,\]
holds. If one writes \smash{$q={\rm e}^{2\pi\sqrt{-1}\sigma}$}, \smash{$p={\rm e}^{2\pi\sqrt{-1}\tau}$},
with complex $\sigma$, $\tau$, then $g(u)$ is indeed periodic in $u$
with periods $\sigma^{-1}$ and $\tau\sigma^{-1}$
(which can be verified by applying \eqref{p1id} to each of the $2s$
theta functions appearing in $g(u)$).
Keeping this notation for $p$ and $q$, we denote the {\em field
of elliptic functions} over $\C$ of the complex variable $u$,
with the two periods $\sigma^{-1}$ and
$\tau\sigma^{-1}$, by $\E_{q,p}(q^u)$.
(We use the notation $\E_{q,p}(q^u)$ instead of
$\E_{q,p}(u)$ as we work with \textit{multiplicatively} denoted
theta functions.)
More generally, we denote the {\em field of totally elliptic multivariate
functions} over $\C$ of the complex variables $u_1,\dots,u_n$,
in each variable with two equal periods,
$\sigma^{-1}$ and $\tau\sigma^{-1}$, by
$\E_{q,p}(q^{u_1},\dots,q^{u_n})$.
The notion of totally elliptic multivariate functions
was first introduced by Spiridonov~\mbox{\cite{Sp02,Sp11}}.

Recall that an (ordinary) \textit{hypergeometric series} is
a series $\sum_{k\ge 0}c_k$ with $c_0=1$ such that $g(k):=c_{k+1}/c_k$
is a rational function in $k$. Further, a \textit{basic hypergeometric series}
(also called \textit{$q$-hypergeometric series}) is
a series \smash{$\sum_{k\ge 0}c_k$} with $c_0=1$ such that $g(k):=c_{k+1}/c_k$
is a rational function in $q^k$. Similarly, an \textit{elliptic
hypergeometric series} is defined to be a series $\sum_{k\ge 0}c_k$
with~${c_0=1}$ such that $g(k):=c_{k+1}/c_k$ is an \textit{elliptic function}
in $k$ (viewed as a complex variable).
The definition of an elliptic hypergeometric series extends that
of a basic hypergeometric series, assuming that
the rational function in $q^k$ that appears in the definition
of the basic hypergeometric series is a ratio
of two polynomials $\alpha\bigl(q^k\bigr)=\sum_{j=0}^s\alpha_jq^{kj}$ and
$\beta\bigl(q^k\bigr)=\sum_{j=0}^s\beta_jq^{kj}$
of equal degree $s$, with non-vanishing constant terms
such that the (polynomial) balancing condition~${\alpha_s/\alpha_0
=\beta_s/\beta_0}$ holds.

We conclude our brief introduction by explicitly reproducing
Frenkel and Turaev's ${}_{10}V_9$ summation \cite{FT}
(see also \cite[equation~(11.4.1)]{GR}),
an identity which is fundamental to the
theory of elliptic hypergeometric series:
Let $m\in\N_0$ and $a,b,c,d,e,q,p\in\C$ with $|p|<1$.
Then there holds the following identity:
\begin{gather}
\sum_{k=0}^m\frac{\ta\bigl(aq^{2k}\bigr)}{\ta(a)}
\frac{(a,b,c,d,e,q^{-m};q,p)_k}
{\bigl(q,aq/b,aq/c,aq/d,aq/e,aq^{m+1};q,p\bigr)_k}q^k\nonumber\\
\qquad
=\frac{(aq,aq/bc,aq/bd,aq/cd;q,p)_m}
{(aq/b,aq/c,aq/d,aq/bcd;q,p)_m},\label{propfteq}
\end{gather}
where $a^2q^{m+1}=bcde$.
It is easy to see that the series in~\eqref{propfteq}
is indeed an elliptic hypergeometric series.
The convention of referring to the above series as a
${}_{10}V_9$ series follows the arguments
in~\cite{Sp02b} and has become standard
(see also \cite[Chapter~11]{GR}).

The rest of the paper is organized as follows:
In Section~\ref{sec:bmcoeffs}, we introduce the specific
elliptic weights which we use and define corresponding
elliptic binomial and multinomial coefficients.
The elliptic binomial coefficients considered
here are those which we introduced in \cite{S}
in the context of lattice path enumeration.
They also appeared as the normal form coefficients in an
elliptic extension of the binomial theorem, featured in
\cite[Section~4]{S1}. (Different elliptic binomial coefficients
were considered by Rains~\cite[Definition~11]{R06}
and implicitly also by Coskun and Gustafson~\cite{CG},
both in the context of convolutions for families of multivariate
special functions that are recursively defined by vanishing
properties and a branching rule.)
The elliptic multinomial coefficients, defined in the same section,
extend our elliptic binomial coefficients and are new.
What is interesting about our specific elliptic weights in \eqref{wdef}
is that they satisfy a~certain consistency relation in \eqref{eq:estr}
that is reminiscent of the dynamical Yang--Baxter
equation but has a very simple form as it
involves only scalars. Our analysis in the subsequent section
crucially depends on this consistency relation.
In Section~\ref{sec:ellalg}, we introduce an
algebra of elliptic commuting variables for which we
identify a natural basis. We also look at specific
commutation relations useful for normalization
of the elements of the algebra.
We then turn to another main result of the paper in
Section~\ref{sec:ellmthm}, featuring an elliptic extension
of the multinomial theorem, valid as an identity in the
just introduced algebra of elliptic commuting variables.
The new result extends the elliptic binomial theorem from
\cite[Section~4]{S1} to higher rank. In Section~\ref{sec:conv},
we show how our elliptic
multinomial theorem can be used to rederive
Rosengren's~\cite[Theorem~5.1]{R04} $\mathsf A_r$
extension of the Frenkel--Turaev summation,
which in the basic case was first obtained by
Milne~\cite[Theorem~6.17]{M88}.
In~a~concluding remark we explain how our algebraic derivation
of the $\mathsf A_r$ Frenkel--Turaev summation admits
a direct combinatorial interpretation in terms of elliptic
weighted lattice paths in the integer lattice $\mathbb Z^r$.

\section{Elliptic weights, elliptic binomial and
 multinomial coefficients}\label{sec:bmcoeffs}

For indeterminates $a$, $b$, complex numbers $q$, $p$ (with $|p|<1$),
and integers $s$, $t$, we define the \textit{small elliptic weights} by
\begin{gather}\label{wdef}
w_{a,b;q,p}(s,t):=
\frac{\ta\bigl(aq^{s+2t},bq^{2s+t-2},aq^{t-s-1}/b\bigr)}
{\ta\bigl(aq^{s+2t-2},bq^{2s+t},aq^{t-s+1}/b\bigr)}q.
\end{gather}
The corresponding \textit{big elliptic weights} are defined by
\begin{gather}\label{Wdef}
W_{a,b;q,p}(s,t):=\prod_{j=1}^tw_{a,b;q,p}(s,j)
=\frac{\ta\bigl(aq^{s+2t},bq^{2s},bq^{2s-1},aq^{1-s}/b,aq^{-s}/b\bigr)}
{\ta\bigl(aq^s,bq^{2s+t},bq^{2s+t-1},aq^{1+t-s}/b,aq^{t-s}/b\bigr)}q^t.
\end{gather}
Clearly, $W_{a,b;q,p}(s,0)=1$, for all $s$.
For $p\to 0$ followed by $a\to 0$ and $b\to 0$
(or $p\to 0$ followed by $b\to \infty$ and $a\to \infty$),
the small elliptic weights $w_{a,b;q,p}(s,t)$ all reduce to $q$
and the big elliptic weights $W_{a,b;q,p}(s,t)$ reduce to $q^t$.
For convenience, we also define the following shifted variant of
a big elliptic weight,
\begin{gather}\label{sWdef}
W_{a,b;q,p}^{(\rho)}(s,t):=W_{aq^{2\rho},bq^{2\rho};q,p}(s,t),
\end{gather}
and further the \textit{big $Q$-weights} by the product
\begin{align}
Q_{a,b;q,p}(\ell,\rho,s,t):={}& \prod_{i=1}^\ell
W_{a,b;q,p}^{(\rho)}(i+s,t)\label{Qweights}\\
={}&\frac{\bigl(aq^{1+2\rho+s+t};q,p\bigr)_\ell \bigl(bq^{1+2\rho+2s};q,p\bigr)_{2\ell}
\bigl(aq^{1-\ell-s}/b,aq^{-\ell-s}/b;q,p\bigr)_\ell}
{\bigl(aq^{1+2\rho+s};q,p\bigr)_\ell \bigl(bq^{1+2\rho+2s+t};q,p\bigr)_{2\ell}
\bigl(aq^{1+t-\ell-s}/b,aq^{t-\ell-s}/b;q,p)_\ell}q^{\ell t}.\nonumber
\end{align}

Assuming $c$ to be an additional indeterminate,
we would like to highlight the following relation satisfied by the
small elliptic weights \eqref{wdef}, for all $s$,
which we refer to as the \textit{elliptic consistency
 relation}
\begin{gather}
w_{aq^2,bq^2;q,p}(s,s) w_{a,c;q,p}(s,s) w_{bq^2,cq^2;q,p}(s,s)=w_{b,c;q,p}(s,s) w_{aq^2,cq^2;q,p}(s,s) w_{a,b;q,p}(s,s).\label{eq:estr}
\end{gather}
This specific equality of simple products
is readily verified using the explicit expression
for the small elliptic weights in \eqref{wdef}.
Equation~\eqref{eq:estr} concerns an equality
involving six weights in total. In three of the six weights
the ``dynamical'' parameters $a$, $b$, $c$ are shifted by $q^2$,
while in the other three weights those parameters are not shifted.
If we view the weights in \eqref{eq:estr} as nodes in a graph,
and connect the shifted weights by edges from left to right
and the similarly do this separately for the unshifted weights,
an overlapping of two chains becomes apparent, see Figure~\ref{fig1}.
This picture helps to remember the special form of~\eqref{eq:estr}.
\begin{figure}[h]\centering

\setlength{\unitlength}{1.4pt}
\thinlines
\begin{picture}(120,40)
\put(20,20){\line(4,-1){40}}
\put(60,10){\line(4,1){40}}
\put(20,10){\line(4,1){17}}
\put(100,10){\line(-4,1){17}}
\put(60,20){\line(4,-1){17}}
\put(60,20){\line(-4,-1){17}}
\multiput(20,20)(40,0){3}{\circle*{1.8}}
\multiput(20,10)(40,0){3}{\circle*{1.8}}
\put(-20,26){$w_{aq^2,bq^2;q,p}(s,s) w_{a,c;q,p}(s,s) w_{bq^2,cq^2;q,p}(s,s)$}
\put(-20,1){\quad$w_{b,c;q,p}(s,s) w_{aq^2,cq^2;q,p}(s,s) w_{a,b;q,p}(s,s)$}
\end{picture}
 \caption{Interlacing of the shifted terms
 in the elliptic consistency relation.}\label{fig1}
\end{figure}

The elliptic consistence relation in \eqref{eq:estr}
guarantees the unique normalization of $X_kX_jX_i$ (for~${1\le i<j<k\le r}$) in the elliptic
algebra that is defined in Definition~\ref{defear},
in particular it is responsible for the equality of the two
expressions obtained on the right-hand sides of
equations~\eqref{eq:diamond1} and \eqref{eq:diamond2}.
Equation~\eqref{eq:estr} is reminiscent
of the dynamical Yang--Baxter equation (a~master equation
in integrable models in statistical mechanics and quantum
field theory, see~\cite{J89}) but has a very simple form,
as it involves only scalars (or $1\times 1$ matrices)
and no operators.
Independently, to the best of our knowledge, it was not
known before that \eqref{eq:estr} has the solution~\eqref{wdef},
not even in the case $p=0$.

For indeterminates $a$, $b$, complex numbers $q$, $p$
(with $|p|<1$), and integers $n$, $k$, we
define the {\em elliptic binomial coefficient} as follows
\begin{gather}\label{ellbc}
\begin{bmatrix}n\\k\end{bmatrix}_{a,b;q,p}:=
\frac{\bigl(q^{1+k},aq^{1+k},bq^{1+k},aq^{1-k}/b;q,p\bigr)_{n-k}}
{\bigl(q,aq,bq^{1+2k},aq/b;q,p\bigr)_{n-k}}.
\end{gather}
This is exactly the expression for $w(\mathcal P((0,0)\to(k,n-k)))$
in \cite[Theorem~2.1]{S}. In \cite{S1}, it was moreover shown that the
elliptic binomial coefficients in~\eqref{ellbc} indeed appear as the
coefficients in a noncommutative elliptic binomial theorem.
Note that the elliptic binomial coefficient
in~\eqref{ellbc} reduces to the usual $q$-binomial coefficient
after taking the limits $p\to 0$, $a\to 0$, and $b\to 0$,
in this order (or after taking the limits in the order $p\to 0$,
$b\to\infty$, and $a\to\infty$).
As pointed out in \cite{S}, the expression in
\eqref{ellbc} is {\em totally elliptic}, i.e.,
elliptic in each of $\log_qa$, $\log_qb$, $n$, and $k$
(viewed as complex parameters), with equal periods of double periodicity.
In particular,
$\left[\begin{smallmatrix}n\\k\end{smallmatrix}\right]_{a,b;q,p}
\in\E_{q,p}\bigl(a,b,q^n,q^k\bigr)$.

It is immediate from the definition of \eqref{ellbc} that, for
integers $n$, $k$, there holds
\begin{subequations}\label{qbinrel}
\begin{gather}
\begin{bmatrix}n\\0\end{bmatrix}_{a,b;q,p}=
\begin{bmatrix}n\\n\end{bmatrix}_{a,b;q,p}=1,
\end{gather}
and
\begin{gather}
\begin{bmatrix}n\\k\end{bmatrix}_{a,b;q,p}=0,\qquad\text{whenever}\
k<0,\ \text{or}\ k> n.
\end{gather}
Furthermore, using the theta addition formula in \eqref{addf}
one can verify the following recursion formula for the
elliptic binomial coefficients:
\begin{gather}
\begin{bmatrix}n+1\\k\end{bmatrix}_{a,b;q,p}=
\begin{bmatrix}n\\k\end{bmatrix}_{a,b;q,p}+
\begin{bmatrix}n\\k-1\end{bmatrix}_{a,b;q,p} W_{a,b;q,p}(k,n+1-k),
\end{gather}
\end{subequations}
for non-negative integers $n$ and $k$.

In the above classical limit,
the relations in \eqref{qbinrel} reduce to
\begin{gather*}
\begin{bmatrix}n\\0\end{bmatrix}_{q}=
\begin{bmatrix}n\\n\end{bmatrix}_{q}=1,
\qquad
\begin{bmatrix}n+1\\k\end{bmatrix}_{q}=
\begin{bmatrix}n\\k\end{bmatrix}_{q}+
\begin{bmatrix}n\\k-1\end{bmatrix}_{q} q^{n+1-k},
\end{gather*}
for positive integers $n$ and $k$ with $n\ge k$, which is
a well-known recursion for the $q$-binomial coefficients.

As was shown in \cite{S1}, the elliptic binomial coefficients
in \eqref{ellbc} can be interpreted as the (area) generating function
for all lattice paths in the integer lattice $\Z^2$ from $(0,0)$ to
$(k,n-k)$ consisting of East and North steps of unit length
where each path is weighted with respect to the product of the
weights of the respective squares covered by the path.
In this interpretation, the weight of the single square with
north-east corner $(s,t)$ is given by $w_{a,b;q,p}(s,t)$,
whereas~$W_{a,b;q,p}(s,t)$ can be regarded as
the weight of the $s$-th column having height $t$.

To prepare the reader for a better understanding of the main result
of this paper, namely the elliptic multinomial theorem in
Section~\ref{sec:ellmthm}, it will be convenient to recall the
author's elliptic binomial theorem from \cite[Theorem~2]{S1}.
We start with the definition of
the algebra of elliptic commuting variables in which the
elliptic binomial coefficients manifestly appear as the
coefficients in a binomial expansion after normal ordering of the
respective variables.

\begin{Definition}\label{defea}
For two nonzero complex numbers $q$ and $p$ with $|p|<1$,
let $\C_{q,p}[X,Y;a,b]$ denote
the associative unital algebra over $\C$,
generated by $X$, $Y$,
satisfying the following three relations:
\begin{align*}
YX=W_{a,b;q,p}(1,1) XY,
\qquad
Xf(a,b)=f\bigl(aq,bq^2\bigr) X,\qquad
Yf(a,b)=f\bigl(aq^2,bq\bigr) Y,
\end{align*}
for all $f\in\E_{q,p}(a,b)$.
\end{Definition}
We refer to the variables $X$, $Y$, $a$, $b$
forming $\C_{q,p}[X,Y;a,b]$
as {\em elliptic commuting} variables.
The algebra $\C_{q,p}[X,Y;a,b]$ reduces to
$\C_{q}[X,Y]$ if one formally lets $p\to 0$, $a\to 0$,
then $b\to 0$, in this order, or lets $p\to 0$, $b\to\infty$,
then $a\to\infty$, in this order, while (having eliminated the nome $p$)
relaxing the condition of ellipticity.
It should be noted that the monomials $X^kY^l$ form a basis for the
algebra $\C_{q,p}[X,Y;a,b]$ as a left module
over $\E_{q,p}(a,b)$,
i.e., any element can be written uniquely as a finite sum
$\sum_{k,l\ge 0} c_{kl}X^kY^l$ with $c_{kl}\in \E_{q,p}(a,b)$
which we call the \textit{normal form} of the element.

The following result from \cite[Theorem~2]{S1}
shows that the normal form of the binomial
${(X+Y)^n}$ is ``nice''; each coefficient to the left of
$X^kY^{n-k}$ completely factorizes as an expression in~$\E_{q,p}(a,b)$.

\begin{Theorem}[binomial theorem for 
 variables in ${\C_{q,p}[X,Y;a,b]}$]\label{ebthm}
 Let $n\in\N_0$.
 Then, as an identity in $\C_{q,p}[X,Y;a,b]$, we have
\begin{gather*}
(X+Y)^n=\sum_{k=0}^n\begin{bmatrix}n\\k\end{bmatrix}_{a,b;q,p}X^kY^{n-k}.
\end{gather*}
\end{Theorem}

In \cite[Corollary~4]{S1}, convolution was applied to this result
(together with comparison of coefficients)
yielding the Frenkel and Turaev ${}_{10}V_9$ summation \cite{FT}
in a form equivalent to \eqref{propfteq} by analytic continuation.

\begin{Remark}\label{remHKKS}
In the recent work \cite[Definition~5.6 and Theorem~5.7]{HKKS}, the
author, in collaboration with Hoshi, Katori, and Koornwinder, defined
a similar but different elliptic commuting algebra
with a corresponding binomial theorem.
\end{Remark}

Before we extend the elliptic binomial coefficients in \eqref{ellbc}
to elliptic multinomial coefficients, we rewrite the elliptic
partial fraction decomposition \eqref{pfdA} in a form that
is suitable for our purpose.
Replacing $r$ by $r+1$ in \eqref{pfdA},
isolating the $(r+1)$-th term of the sum, putting the first~$r$~terms to the other side and dividing both sides of the
equation by the $(r+1)$-th term and using~\eqref{eq:inv},
we obtain the following
form of the type $\mathsf A$ elliptic partial fraction identity
\begin{gather}\label{pfd}
 1=\frac{\prod_{j=1}^r\ta(a_{r+1}/a_j)}
 {\prod_{j=1}^{r+1}\ta(a_{r+1}/b_j)}
 \sum_{i=1}^r\frac{\prod_{j=1}^{r+1}\ta(b_j/a_i)}
 {\prod_{\substack{1\le j\le r+1\\j\neq i}}\ta(a_j/a_i)},
\end{gather}
now subject to the elliptic balancing condition
$a_1\cdots a_{r+1}=b_1\cdots b_{r+1}$.

We are ready to define (for the first time)
elliptic multinomial coefficients.
Let $r>1$ be an integer and $a_1,\dots,a_r\in\C^\times$ be variables.
Further, let $k_1,\dots,k_r$ be integers satisfying $k_1+\cdots+k_r\ge 0$.
Here and throughout, we write $K_i:=\sum_{\nu=1}^i k_\nu$,
for $i=0,\dots,r$, and we will later similarly use the notations
$N_i:=\sum_{\nu=1}^i n_\nu$ and $L_i:=\sum_{\nu=1}^i l_\nu$.
We define the elliptic multinomial coefficients explicitly as
\begin{gather}
\begin{bmatrix}k_1+\cdots+k_r\\k_1,\dots,k_r\end{bmatrix}
_{a_1,\dots,a_r;q,p}\nonumber\\
\qquad:=
\frac{(q;q,p)_{k_1+\cdots+k_r}}{\prod_{i=1}^r(q;q,p)_{k_i}}
\prod_{i=1}^r\frac{\bigl(a_iq^{1+K_r-k_i};q,p\bigr)_{k_i}}
{\bigl(a_iq^{1+2K_{i-1}};q,p\bigr)_{k_i}}
\prod_{1\le i<j\le r}\frac{\bigl(a_iq^{1-k_i}/a_j;q,p\bigr)_{k_j}}
{(a_iq/a_j;q,p)_{k_j}}.\label{ellmc}
\end{gather}
For $r=2$, the elliptic multinomial coefficients
\smash{$\big[\begin{smallmatrix}k_1+k_2\\k_1,k_2\end{smallmatrix}\big]_{a_1,a_2;q,p}$}
reduce to the elliptic binomial coefficients
\smash{$\big[\begin{smallmatrix}k_1+k_2\\k_1\end{smallmatrix}\big]_{a_1,a_2;q,p}$}
(which in general is different from
\smash{$\big[\begin{smallmatrix}k_1+k_2\\k_2\end{smallmatrix}\big]_{a_1,a_2;q,p}$})
given in \eqref{ellbc}. That is, for $r=2$ we have two short notations for the
elliptic multinomial coefficients in \eqref{ellmc},
just as in the familiar ordinary case.

The elliptic multinomial coefficients in \eqref{ellmc} satisfy
\[
 \begin{bmatrix}0\\0,\dots,0\end{bmatrix}_{a_1,\dots,a_r;q,p}
=1,
\]
and (remember that we are assuming $k_1+\cdots+k_r\ge 0$)
\[
 \begin{bmatrix}k_1+\cdots+k_r\\k_1,\dots,k_r\end{bmatrix}_{a_1,\dots,a_r;q,p}
=0,\qquad\text{whenever} \ k_j<0 \ \text{for some} \ j=1,\dots,r,
\]
and for $k_1+\cdots+k_r>0$
the recurrence relation
\begin{gather}
\begin{bmatrix}k_1+\cdots+k_r\\
 k_1,\dots,k_r\end{bmatrix}_{a_1,\dots,a_r;q,p}\nonumber\\
 \qquad=\sum_{i=1}^r
 \begin{bmatrix}k_1+\cdots+k_r-1\\k_1,\dots,k_{i-1},k_i-1,
 k_{i+1},\dots,k_r\end{bmatrix}_{a_1,\dots,a_r;q,p}
 \prod_{j>i}W_{a_i,a_j;q,p}^{(K_{j-1}-k_i)}(k_i,k_j).\label{rec-ellmc}
\end{gather}
The latter is readily established by using the elliptic
partial fraction decomposition \eqref{pfd}.
Indeed, dividing both sides of \eqref{rec-ellmc} by
the elliptic multinomial coefficient on the left-hand side and
replacing the elliptic multinomial coefficients and the shifted big
elliptic weights by their explicit expressions in \eqref{ellmc},
\eqref{sWdef} and \eqref{Wdef}, we obtain, after cancellation
of common factors, \eqref{pfd} with respect to the following simultaneous substitutions:
\begin{alignat*}{3}
 & a_i\mapsto q^{k_i}/a_i\quad\text{for}\ 1\le i\le r,\qquad&&
 a_{r+1}\mapsto q^{k_1+\cdots+k_r},&\\
 & b_i\mapsto 1/a_i\quad\text{for}\ 1\le i\le r,\qquad&&
 b_{r+1}\mapsto q^{2(k_1+\cdots+k_r)}.&
\end{alignat*}
This confirms \eqref{rec-ellmc}.

\section{An algebra of elliptic commuting variables}\label{sec:ellalg}

Recall (from Section~\ref{sec:intro})
that $\E_{q,p}(a_1,\dots,a_r)$
denotes the field of totally
elliptic functions over $\C$, in the complex variables
$\log_qa_i$, $1\le i\le r$, with equal periods $\sigma^{-1}$,
$\tau\sigma^{-1}$ (where \smash{$q={\rm e}^{2\pi\sqrt{-1}\sigma}$}, \smash{$p={\rm e}^{2\pi\sqrt{-1}\tau}$},
$\sigma,\tau\in\C$), of double periodicity.

We shall work in the following algebra.
\begin{Definition}\label{defear}
For $2r$ noncommuting variables $X_1,\dots,X_r$, and $a_1,\dots,a_r$,
where the variables $a_1,\dots,a_r$ commute with each other,
and two nonzero complex numbers $q$, $p$ with $|p|<1$, let $\C_{q,p}[X_1,\dots,X_r;a_1,\dots,a_r]$
denote the associative
unital algebra over $\C$, generated by $X_1,{\dots},X_r$,
satisfying the following relations:
\begin{subequations}\label{defeaeqr}
\begin{gather}
X_jX_i=w_{a_i,a_j;q,p}(1,1) X_iX_j
\qquad\text{for}\ 1\le i<j\le r,\label{elleqr}\\
X_i f(a_1,\dots,a_r)=f\bigl(a_1q^2,\dots,a_{i-1}q^2,a_iq,
a_{i+1}q^2,\dots,a_r q^2\bigr) X_i\qquad\text{for}\ 1\le i\le r,\label{xfr}
\end{gather}
\end{subequations}
for all $f\in\E_{q,p}(a_1,\dots,a_r)$,
and where the elliptic weights
$w_{a_i,a_j;q,p}\in\E_{q,p}(a_1,\dots,a_r)$
are defined in \eqref{wdef}.
\end{Definition}
We refer to the $2r$ variables $X_1,\dots,X_r$, $a_1,\dots,a_r$
forming $\C_{q,p}[X_1,\dots,X_r;a_1,\dots,a_r]$
as {\em elliptic-commuting} variables.\footnote{The algebra
 $\C_{q,p}[X_1,\dots,X_r;a_1,\dots,a_r]$
reduces to the well-known algebra of $q$-commuting
variables, that we may denote by $\C_q[X_1,\dots,X_r]$,
defined by $X_jX_i=qX_iX_j$ for $1\le i<j\le r$,
if one formally lets $p\to 0$ and~${a_1\to 0,\dots,a_r\to 0}$,
in this order, or lets $p\to 0$ and $a_r\to\infty,\dots,a_1\to\infty$,
in this order, while (having eliminated the nome $p$) relaxing the
conditions of ellipticity.}

The following commutation relations, for $1\le i<j\le r$
and $1\le k\le r$, arise as a consequence of
\eqref{xfr} combined with \eqref{wdef}:
\begin{gather*}
 X_i w_{a_i,a_j;q,p}(s,t)=w_{a_i,a_j;q,p}(s+1,t) X_i,\qquad
 X_j w_{a_i,a_j;q,p}(s,t)=w_{a_i,a_j;q,p}(s,t+1) X_j,\\
 X_k w_{a_i,a_j;q,p}(s,t)=w_{a_iq^2,a_jq^2;q,p}(s,t) X_k
 \qquad\text{for}\ k\neq i\ \text{and}\ k\neq j.
\end{gather*}

\begin{Proposition}\label{prop}
The monomials $X_1^{k_1}\cdots X_r^{k_r}$,
for $k_1,\dots,k_r\in\N_0$, form a \textit{basis} for the
algebra~${\C_{q,p}[X_1,\dots,X_r;a_1,\dots,a_r]}$
as a left module over $\E_{q,p}(a_1,\dots,a_r)$.
In other words, any element of the algebra can be written uniquely
as a finite sum
\[
\sum_{k_1,\dots,k_r\ge 0} c_{k_1,\dots,k_r}X_1^{k_1}\cdots X_r^{k_r}
\]
with $c_{k_1,\dots,k_r}\in \E_{q,p}(a_1,\dots,a_r)$ which
we call the \textit{normal form} of the element.
\end{Proposition}
\begin{proof}
It is clear that any element in
$\C_{q,p}[X_1,\dots,X_r;a_1,\dots,a_r]$
can be put into normal form using the relations in \eqref{defeaeqr}.
What still needs to be shown is that the coefficients $c_{k_1,\dots,k_r}$ to the left of each of the
monomials $X_1^{k_1}\cdots X_r^{k_r}$ are well defined in that
they are independent of the order in which the commutation relations
are applied in the normalization procedure.
In other words, since we are working in a multivariate
noncommutative setting, we show that in the associative algebra
$\C_{q,p}[X_1,\dots,X_r;a_1,\dots,a_r]$ a suitable
variant of Bergman's diamond lemma~\cite{B} applies.
By linearity, it is enough to assume that the element consists
of a finite product of the variables~$X_i$ (in some order and
allowing repetitions) and some elements in
$\E_{q,p}(a_1,\dots,a_r)$. Using \eqref{xfr},
the latter elements can all be moved
to the left of the $X_i$ (possibly creating some shifts in the
$a_j$, $1\le j\le r$, by some powers of $q$) resulting in the
product of a unique element in $\E_{q,p}(a_1,\dots,a_r)$ times
a finite product of the variables $X_i$, which without loss
of generality we may assume to be $X=\prod_{s=1}^m X_{i_s}$
(i.e., we assume $m$ to be the degree of the monomial $X$)
which we may refer to as a \textit{word} of $m$ letters.

We now sketch the reduction algorithm of the word
$X=\prod_{s=1}^m X_{i_s}$ to normal form:
If there are any $1\le s<m$ for which two neighboring variables
(\textit{letters})
$X_{i_s}$ and $X_{i_{s+1}}$ are not
in the right order, i.e., when $i_{s}>i_{s+1}$, then we pick
such an $s$ (the choice of $s$ may not be unique) and apply
a reduction of the form \eqref{elleqr} (to switch the order of
$X_{i_s}$ and $X_{i_{s+1}}$) and subsequently use \eqref{xfr}
sufficiently many times to move any
newly created elliptic weights that depend on~${a_1,\dots,a_r}$
to the most left. This step is repeated until the indices
$i_s$, $1\le s\le m$, are in weakly increasing order.
(In the terminology of the diamond lemma, our total monomial
ordering is thus the lexicographic ordering on the letters,
or equivalently, on the indices of the monomials.)
Now, since at each step there may be several $s$ with
$i_{s}>i_{s+1}$, the order of the described reductions is not
unique, leading to possible \textit{ambiguity} in eventually
arriving at an irreducible form
(which is a form where no further reductions can be applied).
According to the diamond lemma, there are two possible ambiguities
when applying reductions:
\textit{overlap ambiguity} and \textit{inclusion ambiguity}.
The special form of the defining commutation relations
\eqref{elleqr} (each term is transformed to another term,
there are no new terms created) shows that in our setting no
inclusion ambiguity can ever arise.
(See~\cite{B} for examples of rings when inclusion ambiguity can arise.)
What still needs to be done is to show that each overlap ambiguity is
\textit{resolvable}. In the quadratic algebra we are considering,
this reduces by the diamond lemma to the following:
We only have to show the any subword of \textit{three} variables
$X=X_kX_jX_i$ with $1\le i<j<k\le r$ is reduction-unique.
See Figure~\ref{fig2} for illustration.

\begin{figure}[ht]\centering
\setlength{\unitlength}{1.2pt}
\thinlines
\begin{picture}(200,70)
\put(00,31){\vector(2,1){45}}
\put(00,22){\vector(2,-1){45}}
\put(152,54){\line(2,-1){9}}
\put(164,48){\line(2,-1){9}}
\put(176,42){\line(2,-1){9}}
\put(188,36){\vector(2,-1){9}}
\put(152,00){\line(2,1){9}}
\put(164,06){\line(2,1){9}}
\put(176,12){\line(2,1){9}}
\put(188,18){\vector(2,1){9}}
\put(-40,23){$X_kX_jX_i$}
\put(48,54){$(w_{a_j,a_k;q,p}(1,1)X_jX_k)X_i$}
\put(48,-08){$X_k(w_{a_i,a_j;q,p}(1,1)X_iX_j)$}
\put(200,23){$C$}
\put(176,45){$s_1$}
\put(163,14){$s_2$}
\put(4,45){$\sigma_{j,k}$}
\put(20,14){$\sigma_{i,j}$}
\end{picture}
\vspace{2mm}

 \caption{Diamond lemma: different reductions of the word $X_kX_jX_i$
 lead to a common expression $C$.}\label{fig2}
\end{figure}

In the figure $\sigma_{i,j}$ refers to the application
of the commutation relation in \eqref{elleqr} with indices $i$ and $j$.
There are two ways to apply a simple commutation relation to
$X=X_kX_jX_i$. In this first step of reduction, we can either
apply $\sigma_{j,k}$ or $\sigma_{i,j}$, leading to two different reductions.
We now have to show that when these reductions are further reduced
(by $s_1$ or $s_2$, both representing ordered sequences of commutation
relations), they lead to a common expression $C$.
We can verify the uniqueness of $C$ directly:
On the one hand, we have
\begin{subequations}
\begin{align}
X_kX_jX_i={}&X_k(X_jX_i)=X_k w_{a_i,a_j;q,p}(1,1) X_iX_j\nonumber\\
 ={}&w_{a_iq^2,a_jq^2;q,p}(1,1) X_kX_iX_j\nonumber\\
 ={}&w_{a_iq^2,a_jq^2;q,p}(1,1)
 w_{a_i,a_k;q,p}(1,1) X_iX_kX_j\nonumber\\
={}&w_{a_iq^2,a_jq^2;q,p}(1,1) w_{a_i,a_k;q,p}(1,1) X_i w_{a_j,a_k;q,p}(1,1) X_jX_k\nonumber\\
={}&w_{a_iq^2,a_jq^2;q,p}(1,1) w_{a_i,a_k;q,p}(1,1) w_{a_jq^2,a_kq^2;q,p}(1,1) X_iX_jX_k.\label{eq:diamond1}
\end{align}
On the other hand, we have
\begin{align}
X_kX_jX_i={}&(X_kX_j)X_i=w_{a_j,a_k;q,p}(1,1) X_jX_kX_i\nonumber\\
 ={}&w_{a_j,a_k;q,p}(1,1) X_j w_{a_i,a_k;q,p}(1,1) X_iX_k\nonumber\\
 ={}&w_{a_j,a_k;q,p}(1,1) w_{a_iq^2,a_kq^2;q,p}(1,1) X_jX_iX_k\nonumber\\
={}&w_{a_j,a_k;q,p}(1,1) w_{a_iq^2,a_kq^2;q,p}(1,1) w_{a_i,a_j;q,p}(1,1) X_iX_jX_k.\label{eq:diamond2}
\end{align}
\end{subequations}
Comparison of the coefficients to the left of
$X_iX_jX_k$ in \eqref{eq:diamond1}
and \eqref{eq:diamond2} gives
\begin{gather*}
w_{a_iq^2,a_jq^2;q,p}(1,1) w_{a_i,a_k;q,p}(1,1) w_{a_jq^2,a_kq^2;q,p}(1,1)\nonumber\\
\qquad=w_{a_j,a_k;q,p}(1,1) w_{a_iq^2,a_kq^2;q,p}(1,1) w_{a_i,a_j;q,p}(1,1),
\end{gather*}
which is indeed true by an instance of the elliptic consistency
relation \eqref{eq:estr}. By the diamond lemma, this establishes
that the reduced normalized form of any word is unique.
Thus the reduced normalized form of any element of the algebra
is unique, as claimed.
\end{proof}

The proof of Proposition~\ref{prop} showed that
when normalizing elements of the elliptic algebra~${\C_{q,p}[X_1,\dots,X_r;a_1,\dots,a_r]}$ and hereby
producing products of elliptic weights,
there are different ways to write those products.
For practical purposes, it will be convenient to define a~preference
between the different possible choices. Subsequently, we shall
prefer the expression obtained in \eqref{eq:diamond2}
to that in \eqref{eq:diamond1},
as it contains less shifts of $q$ in the parameters appearing in the weights.
This choice can repeatedly be applied to larger products of weights to
obtain a~product of weights with smallest possible number of shifts of $q$.

For instance, it follows by application of \eqref{eq:estr} and induction
that for any positive integer $l$,
and indices $1\le i_1<i_2<\cdots<i_l\le r$,
the following $l$-variable commutation relation holds
\[
X_{i_{l}}\cdots X_{i_2}X_{i_1}
=\bigg(\prod_{1\le j<k\le l}w_{a_{i_j}q^{2(k-j-1)},a_{i_k}q^{2(k-j-1)};q,p}(1,1)\bigg)
X_{i_1}X_{i_2}\cdots X_{i_{l}}.
\]

For bringing expressions in
$\C_{q,p}[X_1,\dots,X_r;a_1,\dots,a_r]$
into normal form, the following lemma is useful. (While it was not needed
in the proof of Proposition~\ref{prop}, it will be useful in the proofs
of Theorems~\ref{ellmthm} and \ref{thm:cf-Ar-FT}.)

\begin{Lemma}\label{lem:com}
Let $k_1,\dots,k_r$ and $l_1,\dots,l_r$ be non-negative integers.
The following commutation relation holds as an identity in
$\C_{q,p}[X_1,\dots,X_r;a_1,\dots,a_r]$
\begin{gather*}
X_1^{k_1}\cdots X_r^{k_r}X_1^{l_1}\cdots X_r^{l_r}\\
\qquad
=\bigg(\prod_{1\le i<j\le r}Q_{a_i,a_j;q,p}(l_i,K_{j-1}-k_i+L_{i-1},k_i,k_j)\bigg)
X_1^{k_1+l_1}\cdots X_r^{k_r+l_r},
\end{gather*}
where the big $Q$-weights are defined in \eqref{Qweights}.
\end{Lemma}
\begin{proof}
The identity is readily proved by multiple induction using
\begin{gather*}
X_j^{k}X_i^{l}=Q_{a_i,a_j;q,p}(l,0,0,k) X_i^{l}X_j^k,
\end{gather*}
where $1\le i<j\le r$, for any pair of non-negative integers $k$ and $l$
(which is equivalent to the elliptic specialization of \cite[Lemma~1]{S1}),
combined with repeated application of the commutation rule \eqref{xfr}.
\end{proof}

\section{An elliptic multinomial theorem}\label{sec:ellmthm}
We have the following result.
\begin{Theorem}[elliptic multinomial theorem]\label{ellmthm}
 Let $n\in\N_0$.
 Then the following identity is valid in
$\C_{q,p}[X_1,\dots,X_r;a_1,\dots,a_r]$:
\begin{gather*}
(X_1+\cdots+X_r)^n=\sum_{k_1+\cdots+k_r=n}
\begin{bmatrix}n\\k_1,\dots,k_r\end{bmatrix}_{a_1,\dots,a_r;q,p}
X_1^{k_1}\cdots X_r^{k_r}.
\end{gather*}
\end{Theorem}
\begin{proof}
We proceed by induction on $n$. For $n=0$, the formula is trivial.
Now let $n>0$ ($n$~being fixed) and assume that we have already
shown the formula for all non-negative integers less than $n$.
We have (by separating the last factor, applying induction,
applying a special case of Lemma~\ref{lem:com}, shifting the summation,
and finally combining terms using the recurrence relation~\eqref{rec-ellmc})
\begin{align*}
(X_1+\cdots+X_r)^n={}&(X_1+\cdots+X_r)^{n-1}(X_1+\cdots+X_r)\\
={}&\sum_{k_1+\cdots+k_r=n-1}
\begin{bmatrix}n-1\\k_1,\dots,k_r\end{bmatrix}_{a_1,\dots,a_r;q,p}
X_1^{k_1}\cdots X_r^{k_r}(X_1+\cdots+X_r)\\
={}&\sum_{k_1+\cdots+k_r=n-1}
\sum_{i=1}^r\bigg(
\begin{bmatrix}n-1\\k_1,\dots,k_r\end{bmatrix}_{a_1,\dots,a_r;q,p}
\bigg(\prod_{j>i}W_{a_i,a_j;q,p}^{(K_{j-1}-k_i)}(1+k_i,k_j)\bigg)\\*
&\times
X_1^{k_1}\cdots X_{i-1}^{k_{i-1}}X_i^{k_i+1}
X_{i+1}^{k_{i+1}}\cdots X_r^{k_r}\bigg)\\
={}&\sum_{k_1+\cdots+k_r=n}
\sum_{i=1}^r\bigg(
\begin{bmatrix}n-1\\k_1,\dots,k_{i-1},k_i-1,k_{i+1},\dots,k_r\end{bmatrix}_{a_1,\dots,a_r;q,p}\\*
&\times
\bigg(\prod_{j>i}W_{a_i,a_j;q,p}^{(K_{j-1}-k_i)}(k_i,k_j)\bigg)
X_1^{k_1}\cdots X_r^{k_r}\bigg)\\
={}&\sum_{k_1+\cdots+k_r=n}
\begin{bmatrix}n\\k_1,\dots,k_r\end{bmatrix}_{a_1,\dots,a_r;q,p}
X_1^{k_1}\cdots X_r^{k_r},
\end{align*}
which is what was to be shown.
\end{proof}

\section[Rosengren's A\_r extension of the Frenkel--Turaev summation by convolution]{Rosengren's $\boldsymbol{\mathsf A_r}$ extension of the Frenkel--Turaev\\ summation by convolution}\label{sec:conv}

By convolution, applied to the elliptic multinomial theorem
in Theorem~\ref{ellmthm}, we obtain the following result
which turns out to be equivalent to Rosengren's $\mathsf A_r$ extension of
the Frenkel--Turaev $_{10}V_9$ summation.
\begin{Theorem}\label{thm:cf-Ar-FT}
Let $0\le M\le N$ be two integers, and let $n_1,\dots,n_r\in\N_0$
satisfying $N_r=n_1+\cdots+n_r=N$. Then we have
\begin{gather}
\begin{bmatrix}N\\n_1,\dots,n_r\end{bmatrix}_{a_1,\dots,a_r;q,p}\nonumber\\
\qquad =\sum_{k_1+\cdots+k_r=M}\bigg(
 \begin{bmatrix}M\\k_1,\dots,k_r\end{bmatrix}_{a_1,\dots,a_r;q,p}
\begin{bmatrix}N-M\\n_1-k_1,\dots,n_r-k_r\end{bmatrix}_{a_1q^{2M-k_1},\dots,a_rq^{2M-k_r};q,p}\nonumber\\
\phantom{\qquad=}{}\times\prod_{1\le i<j\le r}Q_{a_i,a_j;q,p}
(n_i-k_i,N_{i-1}+K_{j-1}-K_i,k_i,k_j)\bigg).\label{eq:cf-Ar-FT}
\end{gather}
\end{Theorem}
\begin{proof}
 Working in $\C_{q,p}[X_1,\dots,X_r;a_1,\dots,a_r]$,
 we expand
$(X_1+\cdots+X_r)^N=(X_1+\cdots+X_r)^M(X_1+\cdots+X_r)^{N-M}$
in two different ways and extract coefficients to the left of the monomial~${X_1^{n_1}\cdots X_r^{n_r}}$ where $n_1+\cdots+n_r=N$.
On the left-hand side, the expansion is achieved by a~single application
of Theorem~\ref{ellmthm}, which is simply
\[
(X_1+\cdots+X_r)^N=
\sum_{k_1+\cdots+k_r=N}
\begin{bmatrix}N\\k_1,\dots,k_r\end{bmatrix}_{a_1,\dots,a_r;q,p}
X_1^{k_1}\cdots X_r^{k_r},
\]
whose coefficient of $X_1^{n_1}\cdots X_r^{n_r}$ is clearly
\smash{$\left[\begin{smallmatrix}N\\
 n_1,\dots,n_r\end{smallmatrix}\right]_{a_1,\dots,a_r;q,p}$}.
On the right-hand side, we apply~Theorem~\ref{ellmthm} twice
and bring the expression into normal form by multiple applications
of~\eqref{xfr} (to the second elliptic multinomial coefficient)
and finally apply Lemma~\ref{lem:com} (to bring the product of
two monomials into normal form). The details are as follows
\begin{gather*}
(X_1+\cdots+X_r)^M(X_1+\cdots+X_r)^{N-M}\\
\qquad=\sum_{k_1+\cdots+k_r=M}
\begin{bmatrix}M\\k_1,\dots,k_r\end{bmatrix}_{a_1,\dots,a_r;q,p}
X_1^{k_1}\cdots X_r^{k_r}\\
\phantom{\qquad=}{}\times\sum_{l_1+\cdots+l_r=N-M}
\begin{bmatrix}N-M\\l_1,\dots,l_r\end{bmatrix}_{a_1,\dots,a_r;q,p}
X_1^{l_1}\cdots X_r^{l_r}\\
\qquad=\sum_{\substack{k_1+\cdots+k_r=M\\l_1+\cdots+l_r=N-M}}
\biggl(\begin{bmatrix}M\\k_1,\dots,k_r\end{bmatrix}_{a_1,\dots,a_r;q,p}\\
\phantom{\qquad=}{}\times
\begin{bmatrix}N-M\\l_1,\dots,l_r\end{bmatrix}_{a_1q^{2M-k_1},\dots,a_rq^{2M-k_r};q,p}
X_1^{k_1}\cdots X_r^{k_r}X_1^{l_1}\cdots X_r^{l_r}\biggr)\\
\qquad=\sum_{\substack{k_1+\cdots+k_r=M\\l_1+\cdots+l_r=N-M}}\bigg(
\begin{bmatrix}M\\k_1,\dots,k_r\end{bmatrix}_{a_1,\dots,a_r;q,p}
\begin{bmatrix}N-M\\l_1,\dots,l_r\end{bmatrix}_{a_1q^{2M-k_1},\dots,a_rq^{2M-k_r};q,p}\\
\phantom{\qquad=}{}\times
 \bigg(\prod_{1\le i<j\le r}Q_{a_i,a_j;q,p}(l_i,K_{j-1}-k_i+L_{i-1},k_i,k_j)\bigg)
X_1^{k_1+l_1}\cdots X_r^{k_r+l_r}\bigg).
\end{gather*}
Taking coefficients to the left of $x_1^{n_1}\cdots x_r^{n_r}$ evidently
gives the right-hand side of \eqref{eq:cf-Ar-FT}.
\end{proof}

The convolution identity in Theorem~\ref{thm:cf-Ar-FT} can be regarded
as the combinatorial form of the~$\mathsf A_r$ Frenkel--Turaev
summation
\begin{gather}
\frac{(b/a_1,\dots,b/a_{r+1};q,p)_M}{(q,bz_1,\dots,bz_r;q,p)_M}\nonumber\\
\qquad=\sum_{k_1+\cdots+k_r=M}\prod_{1\le i<j\le r}
\frac{q^{k_i} \ta\bigl(z_jq^{k_j-k_i}/z_i\bigr)}{\ta(z_j/z_i)}
\prod_{i=1}^r\frac{\prod_{j=1}^{r+1}(a_jz_i;q,p)_{k_i}}
{(bz_i;q,p)_{k_i}\prod_{j=1}^r(z_iq/z_j;q,p)_{k_i}}.\label{eq:Ar-FT}
\end{gather}
This identity
was first obtained by Rosengren in \cite[Theorem~5.1]{R04},
see also \cite[equation~(11.7.8)]{GR}.
The $r=2$ case of the identity in \eqref{eq:Ar-FT} is the single-sum
Frenkel--Turaev summation in \eqref{propfteq}.
The $p\to 0$ case of the summation in \eqref{eq:Ar-FT}
was derived earlier by Milne~\cite[Theorem~6.17]{M88}.
Now, \eqref{eq:Ar-FT} contains \eqref{eq:cf-Ar-FT} as a special case:
In \eqref{eq:Ar-FT}, perform the following simultaneous substitutions
\begin{alignat*}{3}
& a_i\mapsto q^{-n_i}\quad\text{for}\ 1\le i\le r,\qquad&&
 a_{r+1}\mapsto q^{-M},&\\
& z_i\mapsto 1/a_i\quad\text{for}\ 1\le i\le r,\qquad&&
 b\mapsto q^{-M-N}.&
\end{alignat*}
These substitutions yield~\eqref{eq:cf-Ar-FT} (after some rewriting).
On the contrary, after rewriting the elliptic multinomial coefficients
and weights in~\eqref{eq:cf-Ar-FT} explicitly in terms of products
of theta-shifted factorials,
the restriction that $n_1,\dots,n_r$ are non-negative integer parameters
can be removed by repeated analytic continuation.
This means that~\eqref{eq:cf-Ar-FT} is actually equivalent to~\eqref{eq:Ar-FT}.

\begin{Remark}
While the above derivation of \eqref{eq:cf-Ar-FT}
involved elliptic commuting variables and algebraic manipulations,
it is not difficult to give combinatorial interpretations of the
respective algebraic expressions
in terms of weighted lattice paths in the $r$-dimensional integer lattice
$\mathbb Z^r$. The multinomial $(X_1+\cdots+X_r)^N$
can be interpreted as the generating function for lattice paths
starting at the origin and consisting of $N$ unit steps
where the $i$th of the $r$ different unit steps increases
the $i$th coordinate in $\mathbb Z^r$ by one
while not changing the other coordinates.
In other words, starting at the origin,
after $N$ steps, the path reaches a point in the intersection of
$\mathbb Z^r$ with the hyperplane $z_1+\dots+z_r=N$.
In this interpretation, for any $r$-tuple of non-negative integers
$(k_1,\dots,k_r)$ whose $i$th component is positive,
the weight of the unit step
\[
(k_1,\dots,k_{i-1},k_i-1,k_{i+1},\dots,k_r)\to
(k_1,\dots,k_r)
\]
is defined to be
\begin{gather}\label{eq:weightsr}
 \prod_{i<j\le r}W_{a_i,a_j;q,p}^{(K_{j-1}-k_i)}(k_i,k_j),
\end{gather}
for any $i=1,\dots,r$, in accordance with the recurrence relation of the
elliptic multinomial coefficients in \eqref{rec-ellmc}.
Assuming the weight of a lattice path in $\mathbb Z^r$
to be the product of the weights (which are all of the form
\eqref{eq:weightsr}) of the unit steps it is composed of,
the weighted generating function of the family of all lattice paths
that start at the origin $(0,\dots,0)$ and, after~${N=n_1+\cdots+n_r}$ unit steps, end in $(n_1\dots,n_r)$,
is the elliptic multinonomial coefficient
\[
\begin{bmatrix}N\\n_1,\dots,n_r\end{bmatrix}_{a_1,\dots,a_r;q,p}.
\]
In this lattice path interpretation the convolution in
Theorem~\ref{thm:cf-Ar-FT} then concerns the generating function
of paths that start at the origin $(0,\dots,0)$ and, after
$N=n_1+\cdots+n_r$ unit steps, end exactly in $(n_1\dots,n_r)$
but is refined according to where, after $M$ steps (for fixed $M$
satisfying~${0\le M\le N}$), the path crosses the hyperplane
$z_1+\dots+z_r=M$.

Our derivation of \eqref{eq:cf-Ar-FT}
by convolution (which as we just explained, can be interpreted
in terms of a weighted counting of lattice paths) appears to
constitute the first combinatorial proof of Rosengren's
$\mathsf A_r$ extension of the Frenkel--Turaev summation,
an identity that is of fundamental importance in the theory of
elliptic hypergeometric series associated with root
systems (cf.\ \cite{RW}).
\end{Remark}

\begin{Remark}
 A natural question concerns the possible wider application of the
 methods developed in this paper.
Concretely, it would be interesting to find a higher rank extension of
 the elliptic binomial theorem from
 \cite[Definition~5.6 and Theorem~5.7]{HKKS} that was mentioned
 in Remark~\ref{remHKKS} and to derive a corresponding
 multivariate Frenkel--Turaev summation by convolution, in the same
 way as Theorem~\ref{thm:cf-Ar-FT} was derived in this section.
 We find this an interesting open problem worthwhile to pursue.
 It is not clear whether the such obtained identity would be equivalent
 to Rosengren's $\mathsf A_r$ extension of the Frenkel--Turaev
 summation or whether
 it would be of a different type such as one of the multivariate
 Frenkel--Turaev summations listed in \cite{RW}.
 One could also ask whether the elliptic binomial coefficients from
 \cite{CG} and \cite{R06} can be identified as the normal form coefficients
 in a suitably defined algebra of elliptic commuting variables
 (with developments parallel to those of this paper).

Furthermore, one can simply ask whether any of the other
 multivariate Frenkel--Turaev summations in \cite{RW} admit similar
 algebraic or combinatorial interpretations as \eqref{eq:Ar-FT} does.
 We~currently have no idea whether this is possible and would find
 this question rather difficult (or challenging) to answer affirmatively.
 Specifically, when considering lattice paths in the $\mathsf C_r$ case
 (see \cite[Theorem~7.1]{R04} for Rosengren's $\mathsf C_r$ extension of
 the Frenkel--Turaev sum), it may well be that instead of
 considering positively directly lattice paths in $\Z^r$ that are bounded
 by a~hyperplane, one would have to allow paths that
 move in the direction (positive or negative) of any axis.
It is likely that one would then work in an associative
algebra that contains noncommuting variables
$C_1,\dots,C_r$, and $X_1,X_1^{-1},\dots,X_r,X_r^{-1}$,
and seek an expansion of the product~${\prod_{i=1}^r\bigl(X_i+C_iX_i^{-1}\bigr)^{n_i}}$ in terms of normalized monomials.
 At the moment this is all speculative.
 Further research is needed to determine whether the methods
 developed in this paper can indeed be applied in the setting
 of root systems other than $\mathsf{A}_r$.
\end{Remark}

\subsection*{Acknowledgements}
The author's research was partly supported by
FWF Austrian Science Fund grant
 \href{https://doi.org/10.55776/P32305}{doi:10.55776/\allowbreak P32305}.

\pdfbookmark[1]{References}{ref}
\LastPageEnding

\end{document}